\documentclass[12pt, reqno]{amsart}
\usepackage{amsmath, amsthm, amscd, amsfonts, amssymb, graphicx, color}
\usepackage[bookmarksnumbered, colorlinks, plainpages]{hyperref}

\hypersetup{colorlinks=true,linkcolor=red, anchorcolor=green, citecolor=cyan, urlcolor=red, filecolor=magenta, pdftoolbar=true}
\newtheorem{theorem}{Theorem}[section]

\newtheorem{proposition}[theorem]{Proposition}
\newtheorem{corollary}[theorem]{Corollary}
\theoremstyle{definition}
\newtheorem{definition}[theorem]{Definition}

\theoremstyle{remark}
\newtheorem{remark}[theorem]{Remark}
\numberwithin{equation}{section}

\begin{document}

\setcounter{page}{1}

\title[Representation theorems for operators on Free Banach spaces]{Representation theorems for operators on Free Banach spaces of countable type}

\author[Aguayo, J., Nova. M, \MakeLowercase{and} Ojeda J.]{J. Aguayo$^1$, M. Nova$^2$ \MakeLowercase{and} J. Ojeda$^1$$^{*}$}

\address{$^{1}$Departamento de Matem\'{a}tica, Facultad de Ciencias F\'{\i}sicas y Matem\'{a}ticas, Universidad de Concepci\'{o}n, Casilla 160-C,
Concepci\'{o}n, Chile.}
\email{\textcolor[rgb]{0.00,0.00,0.84}{jaguayo@udec.cl; jacqojeda@udec.cl}}

\address{$^{2}$Departamento de Matem\'{a}tica y F\'{\i}sica Aplicadas, Facultad de Ingenier\'{\i}a, Universidad Cat\'{o}lica de la Sant\'{\i}sima
Concepci\'{o}n, Casilla 297, Concepci\'{o}n, Chile.}
\email{\textcolor[rgb]{0.00,0.00,0.84}{mnova@ucsc.cl}}



\date{Received: xxxxxx; Revised: yyyyyy; Accepted: zzzzzz.
\newline \indent $^{*}$Corresponding author}

\thanks{This work was partially supported by Proyecto VRID N$^{\circ}$
214.014.038-1.0IN}

\begin{abstract}
This work will be centered in commutative Banach subalgebras of the algebra
of bounded linear operators defined on a Free Banach spaces of countable
type. The main goal of this work wil be to formulate a representation
theorem for these operators through integrals defined by spectral measures
type. In order to get this objective, we will show that, under special
conditions, each one of these algebras is isometrically isomorphic to some
space of continuous functions defined over a compact set. Then, we will
identify such compact developing the Gelfand space theory in the
non-archimedean setting. This fact will allow us to define a measure which
is known as spectral measure. As a second goal, we will formulate a matrix
representation theorem for this class of operators whose entries of these
matrices will be integrals coming from scalar measures.
\end{abstract}

\maketitle

\section{Introduction and notation}

Many researchers have tried to generalize the elemental studies of Banach
algebras from classical case to vectorial structures over non-archimedean
fields. The first big task was to find a result similar to the Gelfand-Mazur
Theorem in this context. But, this theorem failed since every field $\mathbb{%
K}$ with a "non-archimedean valuation" is contained in another field $%
\widetilde{\mathbb{K}}$ whose valuation is an extension of previous one and
both fields are different.

One of the pioneers in the study of non-archimedean Banach algebras of
linear operators and spectral theory in this context has been M. Vishik \cite{VI}, especially in the class of linear operators which admit
compact spectrum. We can also mention another important pioneer, V.
Berkovick \cite{BE}, who made a deep study of this subject on his
survey.

This work will be centered in subalgebras commutative Banach of the algebra
of bounded linear operators defined on a Free Banach spaces of countable
type. The main goal of this work will be to formulate a representation
theorem for these operators through integrals defined by spectral measures
type. In order to get this objective, we will show that, under special
conditions, each one of these algebras is isometrically isomorphic to some
space of continuous functions defined over a compact. Then, we will identify
such compact developing the Gelfand space theory in the non-archimedean
setting. This fact will allow us to define a measure which is known as
spectral measure. As a second goal, we will be to formulate a matrix
representation theorem for this class of operators whose entries of these
matrices will be integrals coming from scalar measures.

Throughout this paper, $\mathbb{K}$ denotes a complete, non-archimedean
valued field and its residue class field is formally real.

In the classical situation we can distinguish two type of normed spaces:
those spaces which are separable and those which are not. If $E$ is a
separable normed space over $\mathbb{K}$, then each one-dimensional subspace
of $E$ is homeomorphic to $\mathbb{K}$, so $\mathbb{K}$ must be separable
too. Nevertheless, we know that there exist non-archimedean fields which are
not separable. Thus, for non-archimedean normed spaces the concept of
separability is meaningless if $\mathbb{K}$ is not separable. However,
linearizing the notion of separability, we obtain a useful generalization of
this concept. A normed space $E$ over a non-archimedean valued field is said
to be of\ countable\ type if it contains a countable subset whose linear
hull is dense in $E$. An example of a normed space of countable type is $%
\left( c_{0},\left\Vert \cdot\right\Vert _{\infty}\right) ,$ where $c_{0}$
is the Banach space of all sequences $x=\left( a_{n}\right) _{n\in\mathbb{N}}
$, $a_{n}\in\mathbb{K}$, for which $\lim_{n\rightarrow\infty}a_{n}=0$ and
its norm is given by $\left\Vert x\right\Vert _{\infty}=\sup\left\{
\left\vert a_{n}\right\vert :n\in\mathbb{N}\right\} .$

A non-archimedean Banach space $E$ is said to be Free Banach space if there
exists a family $\left\{ e_{i}\right\} _{i\in J}$ of non-null vectors of $E$
such that any element $x$ of $E$ can be written in the form of convergent
sum $x=\sum_{i\in J}x_{i}e_{i},\ x_{i}\in\mathbb{K},$ and $\left\Vert
x\right\Vert =\sup_{i\in J}\left\vert x_{i}\right\vert \left\Vert
e_{i}\right\Vert .$ The family $\left\{ e_{i}\right\} _{i\in J}$ is called
orthogonal basis of $E.$ If $s:J\rightarrow\left( 0,\infty\right) ,$ then an
example of Free Banach space is $c_{0}\left( J,\mathbb{K},s\right) ,$ the
collection of all $x=\left( x_{i}\right) _{i\in J}$ such that for any $%
\epsilon>0,$ the set $\left\{ i\in J:\left\vert x_{i}\right\vert s\left(
i\right) >\epsilon \right\} $ is, at most, finite and $\left\Vert
x\right\Vert =\sup_{i\in J}\left\vert x_{i}\right\vert s\left( i\right) $.

We already know that a Free Banach space $E$ is isometrically isomorphic to $%
c_{0}\left( J,\mathbb{K},s\right) ,$ for some $s:J\rightarrow \left(
0,\infty \right) .$ In particular if a Free Banach space is of countable
type, then it is isometrically isomorphic to $c_{0}\left( \mathbb{N},\mathbb{%
K},s\right) ,$ for some $s:\mathbb{N}\rightarrow \left( 0,\infty \right) .$
Note that if $s\left( i\right) \in \left\vert \mathbb{K}\right\vert ,\ $for
each $i\in \mathbb{N},$ then $E$ is isometrically isomorphic to $c_{0}\left( 
\mathbb{N},\mathbb{K}\right) $ (or $c_{0}$ in short). \ For\ more details
concerning Free Banach spaces, we refer the reader to \cite{DI}.

Now, since residual class field of $\mathbb{K}$ is formally real, the
bilinear form%
\begin{equation*}
\left\langle \cdot ,\cdot \right\rangle :c_{0}\times c_{0}\rightarrow 
\mathbb{K};\ \left\langle x,y\right\rangle =\sum_{i=1}^{\infty }x_{i}y_{i}
\end{equation*}%
is an inner product, $\left\Vert \cdot \right\Vert =\sqrt{\left\vert
\left\langle \cdot ,\cdot \right\rangle \right\vert }$ is a non-archimedean
norm in $c_{0}$ and the supremum norm $\left\Vert \cdot \right\Vert _{\infty
}$ coincides with $\left\Vert \cdot \right\Vert $, that is, $\left\Vert
\cdot \right\Vert =\left\Vert \cdot \right\Vert _{\infty }$ (see \cite{NB}). Therefore, to study Free Banach spaces of countable type it is
enough to study the space $c_{0}.$

If $E$ and $F$ are $\mathbb{K}$-normed spaces, then $\mathcal{L}\left(
E,F\right) $ will be the space consisting of all continuous linear maps from 
$E$ into $F.$ If $F=E,$ then $\mathcal{L}\left( E,E\right) =\mathcal{L}%
\left( E\right) .$ For any $T\in\mathcal{L}\left( E,F\right) ,$ $N\left(
T\right) $ will denote its Kernel and $R\left( T\right) $ its range.

A linear operator $T$ from $E$ into $F$ is said to be compact operator if $%
T\left( B_{E}\right) $ is compactoid, where $B_{E}$ =$\left\{ x\in
E:\left\Vert x\right\Vert \leq 1\right\} $ is the unit ball of $E.$ It was
proved in \cite{RO} that $T$ is compact if and only if, for each $%
\epsilon >0,$ there exists a lineal operator of finite-dimensional range $S$
in $\mathcal{L}\left( E,F\right) $ such that $\left\Vert T-S\right\Vert \leq
\epsilon .$

Since $c_{0}$ is not orthomodular, there exist operators in $\mathcal{L}%
\left( c_{0}\right) $ which do not admit adjoint; for example, the linear
operator $T:c_{0}\rightarrow c_{0}$ defined by $T\left( x\right) =\left(
\sum_{i=1}^{\infty }x_{i}\right) e_{1},\ x=\left( x_{i}\right) _{i\in 
\mathbb{N}}\in c_{0}.$ We will denote by $\mathcal{A}_{0}$ the collection of
all elements of $\mathcal{L}\left( c_{0}\right) $ which admit\ adjoint. A
characterization of the elements of $\mathcal{A}_{0}$ (see \cite{AN1}) is the following:%
\begin{equation*}
\mathcal{A}_{0}=\left\{ T\in \mathcal{L}\left( c_{0}\right) :\forall y\in
c_{0},\ \lim_{i\rightarrow \infty }\left\langle Te_{i},y\right\rangle
=0\right\} .
\end{equation*}%
Of course, $\mathcal{A}_{0}$ is a Banach algebra with unit. \ 

We will understand by a normal projection to any projection $%
P:c_{0}\rightarrow c_{0}$ such that $\left\langle x,y\right\rangle =0$ for
each pair $\left( x,y\right) \in N\left( P\right) \times R\left( P\right) .$
An example of normal projection is 
\begin{equation*}
P\left( \cdot\right) =\frac{\left\langle \cdot,y\right\rangle }{\left\langle
y,y\right\rangle }y 
\end{equation*}
for a fix $y\in c_{0}\smallsetminus\left\{ \theta\right\} $.

Now, for each $a=\left( a_{i}\right) _{i\in\mathbb{N}}\in c_{0},$ the linear
operator $M_{a},$ defined by $M_{a}\left( \cdot\right) =\sum_{i=1}^{\infty
}a_{i}\left\langle \cdot,e_{i}\right\rangle e_{i},$ belongs to $\mathcal{A}%
_{0};$ moreover, 
\begin{equation*}
\lim_{n\rightarrow\infty}\left\Vert M_{a}e_{n}\right\Vert
_{\infty}=\lim_{n\rightarrow\infty}\left\Vert
\sum_{i=1}^{\infty}a_{i}\left\langle e_{n},e_{i}\right\rangle
e_{i}\right\Vert _{\infty}=\lim_{n\rightarrow\infty }\left\vert
a_{n}\right\vert =0, 
\end{equation*}
meanwhile, the identity map $Id$ is also an element of $\mathcal{A}_{0}$, but%
\begin{equation*}
\lim_{n\rightarrow\infty}\left\Vert Id\left( e_{n}\right) \right\Vert
_{\infty}=\lim_{n\rightarrow\infty}\left\Vert e_{n}\right\Vert _{\infty}=1. 
\end{equation*}

Let us denote by $\mathcal{A}_{1}$ the collection of all $T\in\mathcal{L}%
\left( c_{0}\right) $ such that $\lim_{n\rightarrow\infty}Te_{n}=\theta,$
i.e., 
\begin{equation*}
\mathcal{A}_{1}=\left\{ T\in\mathcal{L}\left( c_{0}\right) :\lim
_{n\rightarrow\infty}Te_{n}=\theta\right\} . 
\end{equation*}
From the fact that 
\begin{equation*}
\left\vert \left\langle Te_{n},y\right\rangle \right\vert \leq\left\Vert
Te_{n}\right\Vert _{\infty}\left\Vert y\right\Vert _{\infty}, 
\end{equation*}
we have that $\mathcal{A}_{1}\subsetneqq\mathcal{A}_{0}$ since $Id\notin%
\mathcal{A}_{1}.$

If $T,S\in\mathcal{A}_{1},$ then%
\begin{equation*}
\left\langle S,T\right\rangle =\sum_{n=1}^{\infty}\left\langle S\left(
e_{n}\right) ,T\left( e_{n}\right) \right\rangle 
\end{equation*}
is well-defined, is an inner product in $\mathcal{A}_{1}$ and 
\begin{equation*}
\left\Vert T\right\Vert =\sqrt{\left\vert \left\langle T,T\right\rangle
\right\vert }. 
\end{equation*}

By \cite{DI}, we know that each $T\in \mathcal{L}\left(
c_{0}\right) $ can be represented by $T=\sum_{i,j=1}^{\infty }\alpha
_{i,j}e_{j}^{\prime }\otimes e_{i},$ where $\lim_{i\rightarrow \infty
}\alpha _{i,j}=0,\ $for all $j\in \mathbb{N}$. Also, 
\begin{equation*}
\left\Vert T\right\Vert =\sup \left\{ \left\Vert T\left( e_{i}\right)
\right\Vert _{\infty }:i\in \mathbb{N}\right\} =\sup \left\{ \left\vert
\left\langle T\left( e_{i}\right) ,e_{j}\right\rangle \right\vert :i,j\in 
\mathbb{N}\right\} 
\end{equation*}%
and $T$ is compact if and only if \ \ 
\begin{equation*}
\lim_{j\rightarrow \infty }\sup \left\{ \left\vert \alpha _{i,j}\right\vert
:i\in \mathbb{N}\right\} =0.
\end{equation*}%
Now, note that%
\begin{align*}
\left\Vert Te_{n}\right\Vert _{\infty }& =\left\Vert \left(
\sum_{i,j=1}^{\infty }\alpha _{i,j}e_{j}^{\prime }\otimes e_{i}\right)
\left( e_{n}\right) \right\Vert _{\infty }=\left\Vert \sum_{i,j=1}^{\infty
}\alpha _{i,j}e_{j}^{\prime }\left( e_{n}\right) e_{i}\right\Vert _{\infty }
\\
& =\left\Vert \sum_{i=1}^{\infty }\alpha _{i,n}e_{i}\right\Vert _{\infty
}=\sup \left\{ \left\vert \alpha _{i,n}\right\vert :i\in \mathbb{N}\right\} ,
\end{align*}%
thus, 
\begin{equation*}
T\in \mathcal{A}_{1}\Leftrightarrow \left( T\in \mathcal{A}_{0}\text{ and }T%
\text{ is compact}\right) 
\end{equation*}

For the rest of the paper, let us take a fix orthonormal sequence $\left\{
y^{\left( i\right) }\right\} _{i\in\mathbb{N}}$ in $c_{0},$ that is, $%
\left\langle y^{\left( i\right) },y^{\left( j\right) }\right\rangle =0$,$\
i\neq j,$ and $\left\Vert y^{\left( i\right) }\right\Vert _{\infty }=1.$

The next theorem involves normal projections with compact and self-adjoint
operators. The proof can be found in \cite{AN2}.

\begin{theorem}\label{T1}
If the linear operator $T:c_{0}\rightarrow c_{0}$ is compact and
self-adjoint, then there exists an element $\lambda=\left(
\lambda_{i}\right) _{i\in\mathbb{N}}\in c_{0}$ and an orthonormal sequence $%
\left\{ y^{\left( i\right) }\right\} _{i\in\mathbb{N}}$ in $c_{0}$ such that%
\begin{equation*}
T=\sum_{i=1}^{\infty}\lambda_{i}P_{i}, 
\end{equation*}
where 
\begin{equation*}
P_{i}\left( \cdot\right) =\frac{\left\langle \cdot,y^{\left( i\right)
}\right\rangle }{\left\langle y^{\left( i\right) },y^{\left( i\right)
}\right\rangle }y^{\left( i\right) }
\end{equation*}
is the normal projection defined by $y^{\left( i\right) }.$ Moreover, $%
\left\Vert T\right\Vert =\left\Vert \lambda\right\Vert _{\infty}.$
\end{theorem}

\begin{remark}
\ \ \ \ 

\begin{enumerate}
\item This theorem gives us a characterization for compact and self-adjoint
operators. In fact, it is not hard to see that if we take $\lambda=\left(
\lambda_{i}\right) _{i\in\mathbb{N}}\in c_{0}\ $and an orthonormal sequence $%
\left\{ y^{\left( i\right) }\right\} _{i\in\mathbb{N}}$ in $c_{0}$, the
operator 
\begin{equation*}
T_{\lambda}=\sum_{i=1}^{\infty}\lambda_{i}P_{i}, 
\end{equation*}
is compact, self-adjoint and $\left\Vert T_{\lambda}\right\Vert =\left\Vert
\lambda\right\Vert _{\infty}$, where $P_{i}$ is as in the \break Theorem \ref{T1}.

\item The projection family $\left\{ P_{i}:i\in\mathbb{N}\right\} $ is
orthonormal which implies that is orthogonal in the van Rooij's sense$.$ In
fact, for $i\neq j$ 
\begin{equation*}
\left\langle P_{i},P_{j}\right\rangle =\sum_{n=1}^{\infty}\left\langle
P_{i}\left( e_{n}\right) ,P_{j}\left( e_{n}\right) \right\rangle = \left\langle
y^{i},y^{j}\right\rangle \sum_{n=1}^{\infty} \frac{%
y_{n}^{i}}{\left\langle y^{i},y^{i}\right\rangle }\frac{y_{n}^{j}}{%
\left\langle y^{j},y^{j}\right\rangle } = 0 
\end{equation*}
and, for $i\in\mathbb{N}$,  $\left\Vert P_{i}\right\Vert =1$.

\item It is not difficult to prove that $\sum_{s=1}^{N}P_{s}$ and $%
Id-\sum_{r=1}^{M}P_{r}$ are normal projection, when $P_{s}$ and $P_{r}$ are.
\end{enumerate}
\end{remark}

\section{Algebra of operators}

\subsection{A commutative algebra}

From now on, we will consider a fixed orthonormal family $Y=\left\{
y^{\left( i\right) }\right\} _{i\in\mathbb{N}}$ in $c_{0}.$ We will denote
by $\mathfrak{T}_{Y}(c_{0})$ the collection of all compact operators $T_{\mu
},\ \mu\in c_{0},$ where%
\begin{equation*}
T_{\mu}=\sum_{i=1}^{\infty}\mu_{i}P_{i}
\end{equation*}
As we know, the adjoint $T_{\mu}^{\ast}$ of $T_{\mu}$ is itself and $%
\lim_{n\rightarrow\infty}T_{\mu}\left( e_{n}\right) =0.$ On the other hand,
since $Y$ is orthonormal, $T_{\mu}\left( y^{\left( i\right) }\right)
=\mu_{i}y^{\left( i\right) };$ in other words $\mu_{i},\ i\in\mathbb{N}$, is
an eigenvalue of $T_{\mu}.$ Let us denote by $\sigma\left( T_{\mu}\right) $
the set of eigenvalues of $T_{\mu}.$

Now, the collection $\mathfrak{T}_{Y}(c_{0})\,\ $is a linear space with the
operations%
\begin{equation*}
T_{\lambda}+T_{\mu}=T_{\lambda+\mu};\ \ \alpha T_{\lambda}=T_{\alpha\lambda}
\end{equation*}
On the other hand, since $c_{0}$ is a commutative algebra with the operation 
$\lambda\cdot\mu=\left( \lambda_{i}\mu_{i}\right) ,$ we have%
\begin{equation*}
T_{\lambda}\circ T_{\mu}=T_{\lambda\cdot\mu}=T_{\mu}\circ T_{\lambda}. 
\end{equation*}
In order to simply the notation, $T_{\lambda}\circ T_{\mu}$ will be denoted
by $T_{\lambda}T_{\mu}.\ $\newline
With the operations described above, $\mathfrak{T}_{Y}(c_{0})$ becomes a
commutative algebra without unit. Even more, by the fact that $%
T_{\lambda}=T_{\mu}$ implies $\lambda=\mu,$ the map%
\begin{equation*}
\Lambda:c_{0}\rightarrow\mathfrak{T}_{Y}(c_{0});\ \lambda\longmapsto
\Lambda\left( \lambda\right) =T_{\lambda}
\end{equation*}
is an isometric isomorphism of algebras.

As we know, each algebra $E$ without unit can be transformed in an algebra
with unit by considering the collection $E^{+}=\mathbb{K}\oplus E$ provided
with the usual operations and the multiplication operation defined by%
\begin{equation*}
\left( \alpha,\mu\right) \odot\left( \beta,\nu\right) =\left( \alpha
\beta,\alpha\nu+\beta\mu+\mu\cdot\nu\right) . 
\end{equation*}
The unit of this algebra is $\left( 1,\theta\right) ,$ where $\theta$ is the
null vector of $E$. If $E$ is, in particular, a normed space, then so is $%
E^{+}$ and%
\begin{equation*}
\left\Vert \left( \alpha,\mu\right) \right\Vert =\max\left\{ \left\vert
\alpha\right\vert ,\left\Vert \mu\right\Vert _{\infty}\right\} . 
\end{equation*}

It is known that if $E$ is an algebra with power multiplicative norm, that is%
\begin{equation*}
\left\Vert \nu^{n}\right\Vert =\left\Vert \nu\right\Vert ^{n};\ \nu\in E,\
n\in\mathbb{N}. 
\end{equation*}
then the norm on $E^{+}$ is also power multiplicative. As, an example of
algebra with power multiplicative norm is $c_{0}.$

Now, the commutative Banach algebra $\left( \mathfrak{T}_{Y}(c_{0}),+,\cdot,%
\circ,\left\Vert \cdot\right\Vert \right) $ can be transformed, as above, in
a commutative Banach algebra $\left( \mathfrak{T}_{Y}(c_{0})^{+},+,\cdot,%
\odot,\left\Vert \cdot\right\Vert \right) $ with unit. By the fact that $%
c_{0}$ is isometrically isomorphic to $\mathfrak{T}_{Y}(c_{0})$, $\mathfrak{T%
}_{Y}(c_{0})^{+}$ is isometrically isomorphic to $c_{0}{}^{+}.$

We will denote by $\mathcal{S}_{Y}\left( c_{0}\right) $ the collection of
all linear operators $\alpha Id+T_{\lambda},\ $where $\alpha\in\mathbb{K}$, $%
T_{\lambda}\in\mathfrak{T}_{Y}(c_{0})$ and $Id$ is the identity operator on $%
c_{0}.$ $\mathcal{S}_{Y}\left( c_{0}\right) $ is a normed space and if we
add the operation 
\begin{align*}
\left( \alpha_{1}Id+T_{\mu}\right) \left( \alpha_{2}Id+T_{\nu}\right) &
=\alpha_{1}\alpha_{2}Id+\alpha_{1}T_{\nu}+\alpha_{2}T_{\mu}+T_{\mu}T_{\nu} \\
& =\alpha_{1}\alpha_{2}Id+T_{\alpha_{1}\nu+\alpha_{2}\mu+\mu\nu}
\end{align*}
then $\mathcal{S}_{Y}\left( c_{0}\right) $ is converted in a commutative
algebra with unit.

\begin{theorem}
The algebra $\mathcal{S}_{Y}\left( c_{0}\right) $ is isometrically
isomorphic to $\mathfrak{T}_{Y}(c_{0})^{+}.$ As a consequence, $\mathcal{S}%
_{Y}\left( c_{0}\right) $ is a commutative Banach algebra with unit.

\begin{proof}
We define%
\begin{align*}
\mathfrak{T}_{Y}(c_{0})^{+} & \rightarrow\mathcal{S}_{Y}\left( c_{0}\right)
\\
\left( \alpha,T_{\lambda}\right) & \longmapsto\alpha Id+T_{\lambda}
\end{align*}
Since $\alpha T_{\mu}+\beta T_{\lambda}+T_{\lambda}T_{\mu}=T_{\alpha\mu
+\beta\lambda+\mu\lambda},$ the above transformation is a homomorphism of
algebras$.$ Obviously, this homomorphism is onto; hence it is enough to
prove that it is an isometry. We claim that%
\begin{equation*}
\left\Vert \alpha Id+T_{\lambda}\right\Vert =\left\Vert \left( \alpha
,T_{\lambda}\right) \right\Vert 
\end{equation*}
If $\alpha=0$ or $\left\Vert T_{\lambda}\right\Vert =0$ or $\left\Vert
\alpha Id\right\Vert \neq\left\Vert T_{\lambda}\right\Vert ,$ we are done.
We only need to check it when%
\begin{equation*}
\left\vert \alpha\right\vert =\left\Vert \alpha Id\right\Vert =\left\Vert
T_{\lambda}\right\Vert \neq0. 
\end{equation*}
Of course, 
\begin{equation*}
\left\Vert \alpha Id+T_{\lambda}\right\Vert \leq\max\left\{ \left\vert
\alpha\right\vert ,\left\Vert T_{\lambda}\right\Vert \right\} . 
\end{equation*}
Now, by the compactness of $T_{\lambda}$, 
\begin{equation*}
\lim_{n\rightarrow\infty}T_{\lambda}\left( e_{n}\right) =0. 
\end{equation*}
Thus, there exists $N\in\mathbb{N}$ such that 
\begin{equation*}
n\geq N\Rightarrow\left\Vert T_{\lambda}\left( e_{n}\right) \right\Vert
<\left\vert \alpha\right\vert 
\end{equation*}
Therefore,%
\begin{align*}
\left\Vert \alpha Id+T_{\lambda}\right\Vert & =\sup\left\{ \left\Vert \alpha
e_{n}+T_{\lambda}\left( e_{n}\right) \right\Vert :n\in\mathbb{N}\right\} \\
& =\max\left\{ \left\Vert \alpha e_{1}+T_{\lambda}\left( e_{1}\right)
\right\Vert ,\ \left\Vert \alpha e_{2}+T_{\lambda}\left( e_{2}\right)
\right\Vert ,\ldots,\left\Vert \alpha e_{N-1}+T_{\lambda}\left(
e_{N-1}\right) \right\Vert ,\left\vert \alpha\right\vert \right\} \\
& =\left\vert \alpha\right\vert =\max\left\{ \left\vert \alpha\right\vert
,\left\Vert T_{\lambda}\right\Vert \right\} =\left\Vert \left(
\alpha,T_{\lambda}\right) \right\Vert .
\end{align*}
\end{proof}
\end{theorem}

\begin{remark}
Since $c_{0}^{+}$ is isometrically isomorphic to $\mathfrak{T}%
_{Y}(c_{0})^{+},$ the above theorem says that $c_{0}^{+}$ is isometrically
isomorphic to $\mathcal{S}_{Y}\left( c_{0}\right) .$
\end{remark}

We claim that the usual norm in $\mathcal{S}_{Y}\left( c_{0}\right) $ is
power multiplicative, that is,

\begin{proposition}
$\mathcal{S}_{Y}\left( c_{0}\right) $ is an algebra with power
multiplicative norm.

\begin{proof}
It follows from the fact that $c_{0}^{+}$ is isometrically isomorphic to $%
\mathcal{S}_{Y}(c_{0}).$
\end{proof}
\end{proposition}

\begin{definition}
A commutative Banach algebra $\mathcal{A}$ is called a C-algebra if there
exists a locally compact zero-dimensional Hausdorff space $X$ such that $%
\mathcal{A}$ is isometrically isomorphic to $C_{\infty}\left( X\right) ,$
where $C_{\infty}\left( X\right) $ is the space of all continuous functions
from $X$ into $\mathbb{K}$ which vanishes at infinity.
\end{definition}

As we know $\left\{ e_{j}=\left( \delta _{i,j}\right) _{i\in \mathbb{N}%
}:j\in \mathbb{N}\right\} $, where $\delta _{i,j}$ denotes the Kronecker
symbol, is the canonical basis of $c_{0}\,.$ Since%
\begin{equation*}
e_{j}^{2}=e_{j}
\end{equation*}%
and%
\begin{equation*}
\overline{\left\langle \left\{ e_{j}:j\in \mathbb{N}\right\} \right\rangle }%
=c_{0}
\end{equation*}%
we conclude that the collection of all the idempotent elements of $c_{0}$
with norm less than or equal to 1 is dense in $c_{0}.$ As a consequence, $%
c_{0}$ and $c_{0}^{+}$ are C-algebras (see \cite{RO}).

\begin{theorem}
$\mathcal{S}_{Y}\left( c_{0}\right) $ is a C-algebra with unity.

\begin{proof}
It follows from the fact that $c_{0}^{+}$ is isometrically isomorphic to $%
\mathcal{S}_{Y}\left( c_{0}\right) $.
\end{proof}
\end{theorem}

\begin{remark}
We recall that the spectrum of a commutative Banach algebra $\mathfrak{A}$
is the collection $Sp\left( \mathfrak{A}\right) $ of all nonzero algebra
homomorphisms defined from $\mathfrak{A}$ into $\mathbb{K}$, that is, 
\begin{equation*}
Sp\left( \mathfrak{A}\right) =\left\{ \phi:\mathfrak{A}\rightarrow \mathbb{K}%
:\phi\text{ is a nonzero homomorphism}\right\} . 
\end{equation*}
Note that the natural topology on $Sp\left( \mathfrak{A}\right) $ is induced
by the product topology on $\mathbb{K}^{\mathfrak{A}}$ and also for each $%
\phi\in Sp\left( \mathfrak{A}\right) ,$ $\left\Vert \phi\right\Vert \leq1.$
For any $x\in\mathfrak{A},$ we define%
\begin{equation*}
G_{x}:Sp\left( \mathfrak{A}\right) \rightarrow\mathbb{K},\ \ \phi\longmapsto
G_{x}\left( \phi\right) =\phi\left( x\right) 
\end{equation*}
which is clearly continuous and bounded. Let us denote by 
\begin{equation*}
\left\Vert x\right\Vert _{sp}=\sup\left\{ \left\vert \phi\left( x\right)
\right\vert :\phi\in Sp\left( \mathfrak{A}\right) \right\} =\left\Vert
G_{x}\right\Vert _{\infty}
\end{equation*}
the spectral norm of $x.$ Since 
\begin{equation*}
\left\vert \phi\left( x\right) \right\vert \leq\left\Vert \phi\right\Vert
\left\Vert x\right\Vert \leq\left\Vert x\right\Vert , 
\end{equation*}
we have, in general, that 
\begin{equation*}
\left\Vert x\right\Vert _{sp}\leq\left\Vert x\right\Vert . 
\end{equation*}
Finally, let us denote by $R\left( G_{x}\right) $ the range of $G_{x}.$ The
closure of $R\left( G_{x}\right) $ is called spectrum of $x.$
\end{remark}

L. Narici proved the following result (see \cite{RO}):

\begin{proposition}\label{P2}
A commutative Banach algebra $\mathfrak{A}$ with unit is a $C-$algebra if
and only if its spectrum $Sp\left( \mathfrak{A}\right) $ is compact and its
spectral norm $\left\Vert x\right\Vert _{sp}$ is equal to $\left\Vert
x\right\Vert $, for every $x\in\mathfrak{A}.$
\end{proposition}

\begin{remark}
As $\mathcal{S}_{Y}\left( c_{0}\right) $ satisfies the hypothesis of the
Proposition \ref{P2}, we conclude that $Sp\left( \mathcal{S}_{Y}\left(
c_{0}\right) \right) $ is compact and, for every $H\in\mathcal{S}_{Y}\left(
c_{0}\right) ,$ 
\begin{equation*}
\left\Vert H\right\Vert =\sup_{i\in\mathbb{N}}\left\Vert H\left(
e_{i}\right) \right\Vert =\left\Vert H\right\Vert _{sp}. 
\end{equation*}
\end{remark}

Under the conditions that $\mathcal{S}_{Y}\left( c_{0}\right) $ is a $C$%
-algebra commutative with unity and $Sp\left( \mathcal{S}_{Y}\left(
c_{0}\right) \right) $ is compact, we conclude that $\mathcal{S}_{Y}\left(
c_{0}\right) $ is isometrically isomorphic to the space of all continuous $%
\mathbb{K}$-valued functions defined on $Sp\left( \mathcal{S}_{Y}\left(
c_{0}\right) \right) $ provided by the supremum norm, that is, there exists
an isomorphism of algebras 
\begin{equation*}
\Psi:\mathcal{S}_{Y}\left( c_{0}\right) \rightarrow C\left( Sp\left( 
\mathcal{S}_{Y}\left( c_{0}\right) \right) \right) 
\end{equation*}
such that, for all $H\in\mathcal{S}_{Y}\left( c_{0}\right) ,$ $\left\Vert
H\right\Vert =\left\Vert \Psi\left( H\right) \right\Vert _{\infty}.$

Now, since $\mathcal{S}_{Y}\left( c_{0}\right) $ is the closure of the span
of the collection $\left\{ Id,P_{1},P_{2},\ldots\right\} ,$ we can define
the homomorphism of algebra. Let $n\in\mathbb{N};$ 
\begin{equation*}
\phi_{n}:Span\left\{ Id,P_{1},P_{2},\ldots\right\} \rightarrow\mathbb{K}
\end{equation*}
by 
\begin{equation*}
\phi_{n}\left( P_{i}\right) =\left\{ 
\begin{array}{ccc}
1 & if & n=i \\ 
0 & if & n\neq i%
\end{array}
\right. 
\end{equation*}
and $\phi_{0}\left( P_{i}\right) =0,\ $for every $i\in\mathbb{N}.$ Thus, for
any $H=\alpha_{0} Id + \sum_{i=1}^{k}\alpha_{i}P_{i}\in Span\left\{
Id,P_{1},P_{2},\ldots\right\} ,$ we have%
\begin{equation*}
\left\vert \phi_{n}\left( H\right) \right\vert =\left\vert \alpha_{0} Id +
\sum_{i=1}^{k}\alpha_{i}\phi_{n}\left( P_{i}\right) \right\vert \leq
\max\left\{ |\alpha_{0}|, \left\vert \alpha_{1}\phi_{n}\left( P_{1}\right)
\right\vert , ..., \left\vert \alpha_{k}\phi_{n}\left( P_{k}\right)
\right\vert \right\} \leq\left\Vert H\right\Vert , 
\end{equation*}
that is, $\phi_{n}$ is continuous in $Span\left\{
Id,P_{1},P_{2},\ldots\right\} $. From this, $\phi_{n}$ can be uniquely
extended to a continuous algebra homomorphism from $\mathcal{S}_{Y}\left(
c_{0}\right) $ into $\mathbb{K}.$

Note that if $\phi\in Sp\left( \mathcal{S}_{Y}\left( c_{0}\right) \right) ,$
then $\phi\left( P_{i}\right) =\phi_{n}\left( P_{i}\right) $ for some $n\in%
\mathbb{N}\cup\left\{ 0\right\} $. In other words, $Sp\left( \mathcal{S}%
_{Y}\left( c_{0}\right) \right) =\left\{ \phi_{n}:n\in\mathbb{N}\cup\left\{
0\right\} \right\} .$ Therefore, the function $\Gamma:\mathbb{N}\cup\left\{
0\right\} \rightarrow Sp\left( \mathcal{S}_{Y}\left( c_{0}\right) \right) $
defined by $\Gamma\left( n\right) =\phi_{n}$ is bijective.

If we equip $\mathbb{N}\cup\left\{ 0\right\} =\mathbb{N}^{\ast}$ with the
one-point compactification topology of the discrete space $\mathbb{N}$, then 
$\mathbb{N}^{\ast}$ is homeomorphic to $Sp\left( \mathcal{S}_{Y}\left(
c_{0}\right) \right) .$

Since the $\mathcal{S}_{Y}\left( c_{0}\right) \cong C\left( Sp\left( 
\mathcal{S}_{Y}\left( c_{0}\right) \right) \right) $ and the compact space $%
Sp\left( \mathcal{S}_{Y}\left( c_{0}\right) \right) $ is unique up to
homeomorphism, we have that $\mathcal{S}_{Y}\left( c_{0}\right) \cong
C\left( \mathbb{N}^{\ast}\right) .$

Let us identify the isometric isomorphism $\Psi.$ If we replace $Sp\left( 
\mathcal{S}_{Y}\left( c_{0}\right) \right) $ by $\mathbb{N}^{\ast}$ and
considering the map 
\begin{equation*}
G_{T_{\lambda}}\overset{notation}{=}f_{T_{\lambda}}:\mathbb{N}^{\ast
}\rightarrow\mathbb{K};\ n\longmapsto f_{T_{\lambda}}\left( n\right)
=\lambda_{n}\ , 
\end{equation*}
then, for $H=\alpha_{0}Id+T_{\lambda},$ 
\begin{equation*}
G_{H}:\mathbb{N}^{\ast}\rightarrow\mathbb{K};\ n\longmapsto G_{H}\left(
n\right) =\alpha_{0}+\lambda_{n}=\alpha_{0}+f_{T_{\lambda}}\left( n\right) . 
\end{equation*}
Note that, if $\lambda=e_{n},$ then $T_{e_{n}}=P_{n}$ and therefore $%
G_{P_{n}}=\eta_{\left\{ n\right\} },$ where $\eta_{\left\{ n\right\} }$ is
the $\mathbb{K}$-characteristic function.

From this, we can define the transformation 
\begin{equation*}
G:\mathcal{S}_{Y}\left( c_{0}\right) \rightarrow C\left( \mathbb{N}^{\ast
}\right) ;\ H\longmapsto G_{H}. 
\end{equation*}
Clearly, $G$ is an algebra homomorphism and, for any $H\in\mathcal{S}%
_{Y}\left( c_{0}\right) ,$ $\left\Vert H\right\Vert _{sp}=\left\Vert
G_{H}\right\Vert $ and $\left\Vert H\right\Vert _{sp}=\left\Vert
H\right\Vert ;$ in other words, $G$ is an isometry. This transformation is
very-well known as the Gelfand transformation.

As a consequence, we have the following theorem which is analogous to the
Gelfand-Naimark classical theorem in the non-archimedean setting:

\begin{theorem}
$\mathcal{S}_{Y}\left( c_{0}\right) $ is isometrically isomorphic to $%
C\left( \mathbb{N}^{\ast}\right) $ through $G.$
\end{theorem}

\subsection{Spectral measure}\label{S2.2}

Let $\Omega(\mathbb{N}^{\ast})$ be the Boolean ring of all clopen subsets of 
$\mathbb{N}^{\ast}$. The elements of this ring are classified in two classes
of subcollections: the first one contains finite subsets of $\mathbb{N}$
that we will call them of type 1 and the second one contains the complement
of finite subsets in $\mathbb{N}^{\ast}$ that we will call them of type 2.

For a $C\subset\mathbb{N}^{\ast},$ $\eta_{C}$ denotes the $\mathbb{K}$%
-characteristic function of $C.$ If $C_{1},C_{2}\subset\mathbb{N}^{\ast},$
then%
\begin{align*}
\eta_{C_{1}}\cdot\eta_{C_{2}} & =\eta_{C_{1}\cap C_{2}};\ \
\eta_{C}^{2}=\eta_{C} \\
\eta_{C_{1}}+\eta_{C_{2}} & =\eta_{C_{1}\cup C_{2}},\text{ if }C_{1}\cap
C_{2}=\varnothing.
\end{align*}
Of course, $\eta_{C}$ is continuous if and only if $C\in\Omega\left( \mathbb{%
N}^{\ast}\right) .$

Now, let us take $f\in C\left( \mathbb{N}^{\ast}\right) $ and $\epsilon>0.$
Since the subspace generated by $\left\{ \eta_{\left\{ n\right\} }:n\in%
\mathbb{N}\right\} \cup\left\{ \eta_{\mathbb{N}^{\ast}\setminus\left\{
n_{1},n_{2},\ldots,n_{k}\right\} }:\left\{ n_{1},n_{2},\ldots,n_{k}\right\}
\subset\mathbb{N}\right\} $ is $\left\Vert \cdot\right\Vert _{\infty}-dense$
on $C\left( \mathbb{N}^{\ast}\right) ,$ there exists finite collection $%
\left\{ \alpha_{0},\alpha_{1},\cdots,\alpha_{n}\right\} \subset\mathbb{K}$
such that 
\begin{equation*}
\left\Vert f-\left( \alpha_{0}\eta_{\mathbb{N}^{\ast}\setminus\left\{
n_{1},n_{2},\ldots,n_{k}\right\} }+\sum_{s=1}^{k}\alpha_{s}\eta_{\left\{
n_{s}\right\} }\right) \right\Vert _{\infty}=\left\Vert f-\left( \alpha
_{0}\eta_{\mathbb{N}^{\ast}}+\sum_{s=1}^{k}\lambda_{s}\eta_{\left\{
n_{s}\right\} }\right) \right\Vert _{\infty}<\epsilon. 
\end{equation*}
Let us denote by $\Lambda$ the inverse transformation of $G.$ By the
isometry condition of $\Lambda,$ we get that 
\begin{equation*}
\left\Vert \Lambda f-\Lambda\left( \alpha_{0}\eta_{\mathbb{N}%
^{\ast}}+\sum_{s=1}^{k}\lambda_{s}\eta_{\left\{ n_{s}\right\} }\right)
\right\Vert =\left\Vert \Lambda f-\left(
\alpha_{0}Id+\sum_{s=1}^{k}\lambda_{s}P_{n_{s}}\right) \right\Vert
<\epsilon. 
\end{equation*}

Now, we define the set-function $m:\Omega(\mathbb{N}^{\ast})\rightarrow 
\mathcal{S}_{Y}\left( c_{0}\right) $ by $m(C)=\Lambda(\eta_{C}).$ This
set-function is finitely additive and satisfies
\begin{equation}\label{2.1}
m(C)=\left\{ 
\begin{array}{ccc}
\theta & if & C=\varnothing \\ 
Id & if & C=\mathbb{N}^{\ast} \\ 
\sum_{i=1}^{k}P_{n_{i}} & if & C=\{n_{1},...,n_{k}\}\text{ } \\ 
Id-\sum_{i=1}^{k}P_{n_{i}} & if & C=\mathbb{N}^{\ast}\setminus%
\{n_{1},...,n_{k}\}\,%
\end{array}
\right.
\end{equation}

Also, by the fact that $\Vert m(D)\Vert=1,\ $if $D\in\Omega(\mathbb{N}^{\ast
})\setminus\left\{ \emptyset\right\} ,$ the set 
\begin{equation}\label{2.2}
\left\{ m\left( B\right) :B\in\Omega\left( \mathbb{N}^{\ast}\right) ,\
B\subset C\right\}
\end{equation}

is bounded for any $C\in\Omega(\mathbb{N}^{\ast})\setminus\left\{
\emptyset\right\} .$

On the other hand, if $\left\{ C_{\mu}\right\} _{\mu\in\Gamma}$ is shrinking
on $\Omega\left( \mathbb{N}^{\ast}\right) $ and $\cap_{\mu\in\Gamma}C_{\mu
}=\varnothing,$ then there exists $\mu_{0}\in\Gamma$ such that for $\mu\geq$ 
$\mu_{0},$ $C_{\mu}=\varnothing$ and then 

\begin{equation}\label{2.3}
\lim_{\mu\in\Gamma}m\left( C_{\mu}\right) =0.
\end{equation}

By \eqref{2.1},  \eqref{2.2} and \eqref{2.3}, $m$
is a vector measure in the Katsaras's sense (see \cite{KA}).

Let us take a $C\in\Omega(\mathbb{N}^{\ast}),\ C\neq\varnothing$ and denote
by $D_{C}$\textbf{\ }the collection of all \break $\alpha=\left\{
C_{1},C_{2},\ldots,C_{n};x_{1},x_{2},\ldots,x_{n}\right\} ,$ where $\left\{
C_{k}:k=1,\ldots,n\right\} $ is a clopen partition of $C$ and $x_{k}\in
C_{k}.$ We define a partial order in $\mathcal{D}_{C}$ by $%
\alpha_{1}\geq\alpha_{2}$ if and only if the clopen partition of $C$ in $%
\alpha_{1}$ is a refinement of the clopen partition of $C$ in $\alpha_{2}.$
Thus, $\left( \mathcal{D}_{C},\geq\right) $ is a directed set.

Now, for $f\in C\left( \mathbb{N}^{\ast}\right)$, $C\in\Omega (\mathbb{N}%
^{\ast})$ and $\alpha=\left\{
C_{1},C_{2},\ldots,C_{n};x_{1},x_{2},\ldots,x_{n}\right\} \in\mathcal{D}_{C},
$ we define%
\begin{equation*}
\omega_{\alpha}(f,m,C)=\sum_{k=1}^{n}f\left( x_{k}\right) \Lambda
(\eta_{C_{k}})=\Lambda\left( \sum_{k=1}^{n}f\left( x_{k}\right) \eta
_{C_{k}}\right) . 
\end{equation*}
On the other hand, since $f\eta_{C}$ can be reached by a net 
\begin{equation*}
\left\{ \sum_{k=1}^{n}f\left( x_{k}\right) \eta_{C_{k}}\right\} _{\alpha\in%
\mathcal{D}_{C}}
\end{equation*}
in $\left( C\left( \mathbb{N}^{\ast}\right) ,\left\Vert \cdot\right\Vert
_{\infty}\right) $ and $\Lambda$ is continuous, $\lim_{\alpha\in \mathcal{D}%
_{C}}\omega_{\alpha}(f,m,C)$ exists in $\mathcal{S}_{Y}.$ Therefore, the
operator $\Lambda\left( f\eta_{C}\right) $ can be interpreted as an integral
as follows 
\begin{equation*}
\Lambda\left( f\eta_{C}\right) =\int_{\mathbb{N}^{\ast}}f\eta_{C}dm=%
\int_{C}fdm=\lim_{\alpha\in\mathcal{D}_{C}}\omega_{\alpha}(f,m,C). 
\end{equation*}
In particular, 
\begin{equation*}
\Lambda\left( \eta_{C}\right) =\int_{\mathbb{N}^{\ast}}\eta_{C}dm=m\left(
C\right) =\left\{ 
\begin{array}{ccc}
\sum_{i=1}^{s}P_{n_{s}} & if & C=\{n_{1},...,n_{s}\} \\ 
Id-\sum_{i=1}^{k}P_{n_{i}} & if & C=\mathbb{N}^{\ast}\setminus%
\{n_{1},...,n_{k}\}%
\end{array}
\right. . 
\end{equation*}

\begin{theorem}
Each operator in $\mathcal{S}_{Y}$ is represented as an integral defined by
the projection-valued measure.
\end{theorem}

\subsection{Scalar measures and a matrix representation of an operator}\label{S2.3}

In this subsection we will show that each operator of the $C$-algebra $%
\mathcal{S}_{Y}$ can be associated to an infinite matrix where each entry is
an integral defined by a scalar measures.

Let $x,y\in c_{0}.$ We define $m_{x,y}:\Omega(\mathbb{N}^{\ast})\rightarrow 
\mathbb{K}$ by%
\begin{equation*}
m_{x,y}(C)=\left\langle m(C)x,y\right\rangle . 
\end{equation*}
Clearly, $m_{x,y}$ is a scalar measure in the van Rooij's sense and,
following the same arguments given in the previous subsection, 
\begin{equation*}
\Lambda_{x,y}\left( f\right) =\left\langle \Lambda\left( f\right)
x,y\right\rangle 
\end{equation*}
can be interpreted as an integral, say%
\begin{equation*}
\int_{\mathbb{N}^{\ast}}fdm_{x,y}=\left\langle \left( \int_{\mathbb{N}^{\ast
}}fdm\right) \left( x\right) ,y\right\rangle . 
\end{equation*}

In fact, as in the previous subsection 
\begin{equation*}
\omega_{\alpha}(f,m_{x,y},\mathbb{N}^{\ast})=\left\langle \omega_{\alpha
}(f,m,\mathbb{N}^{\ast})\left( x\right) ,y\right\rangle
=\sum_{k=1}^{n}f\left( x_{k}\right) m_{x,y}\left( \mathbb{N}^{\ast}\right) 
\end{equation*}
and then we can denote by 
\begin{equation*}
\int_{\mathbb{N}^{\ast}}f\ dm_{x,y}=\lim_{\alpha\in\mathcal{D}_{\mathbb{N}%
^{\ast}}}\omega_{\alpha}(f,m_{x,y},\mathbb{N}^{\ast})=\lim_{\alpha \in%
\mathcal{D}_{\mathbb{N}^{\ast}}}\left\langle \omega_{\alpha}(f,m,\mathbb{N}%
^{\ast})\left( x\right) ,y\right\rangle =\left\langle H \left( x\right) ,y\right\rangle . 
\end{equation*}
where $H=\int_{\mathbb{N}^{\ast}}f\ dm\in\mathcal{S}_{Y}\left( c_{0}\right) .
$\newline
A particular case is when $x=e_{i}$ and $y=e_{j}.$ In such a case, the
measure $m_{e_{i},e_{j}}$ will be denoted by $m_{ij}.$

Note that, if $C\in\Omega(\mathbb{N}^{\ast})\setminus\left\{ \emptyset
\right\} ,$ then 
\begin{equation*}
\sup_{i,j\in\mathbb{N}}|m_{ij}(C)|=\sup_{i,j\in\mathbb{N}}|\left\langle
m(C)e_{i},e_{j}\right\rangle |=\Vert m(C)\Vert=1=\Vert\eta_{C}\Vert\ \ \ 
\end{equation*}
We define the linear functional%
\begin{equation*}
\Lambda_{ij}:C\left( \mathbb{N}^{\ast}\right) \rightarrow\mathbb{K}\text{;\
\ }f\rightarrow\Lambda_{ij}\left( f\right) =\int_{\mathbb{N}^{\ast }}f\
dm_{ij}=\left\langle H(e_{i}),e_{j}\right\rangle 
\end{equation*}
where $H=\int_{\mathbb{N}^{\ast}}f\ dm\in\mathcal{S}_{Y}\left( c_{0}\right) .
$\newline
Moreover, 
\begin{equation*}
\sup_{i,j\in\mathbb{N}}\left\vert \Lambda_{ij}(f)\right\vert =\sup _{i,j\in%
\mathbb{N}}\left\vert \left\langle H(e_{i}),e_{j}\right\rangle \right\vert
=\Vert H\Vert=\Vert f\Vert 
\end{equation*}

Let us denote by ${\mathcal{M}}$ the space of all infinite matrices of the
form $\left( \Lambda_{ij}(f)\right) _{i,j\in\mathbb{N}},$ i.e., 
\begin{equation*}
{\mathcal{M}}=\left\{ A\left( f\right) =\left( \Lambda_{ij}(f)\right)
_{i,j\in\mathbb{N}}:f\in C\left( \mathbb{N}^{\ast}\right) \right\} 
\end{equation*}
Clearly, ${\mathcal{M}}$ is a vector space over $\mathbb{K}$ and the
function $\left\Vert \cdot\right\Vert _{{\mathcal{M}}}:{\mathcal{M}}%
\rightarrow \mathbb{R}$ defined by $\left\Vert A\left( f\right) \right\Vert
_{{\mathcal{M}}}=\sup_{i,j\in\mathbb{N}}\left\vert
\Lambda_{ij}(f)\right\vert $ is a non-archimedean norm. On the other hand,
if $f=\eta_{\{k\}}$ or $\eta_{\mathbb{N}^{\ast}\setminus\left\{
k_{1},k_{2},\ldots,k_{n}\right\} }$, then 
\begin{align*}
A(\eta_{\{k\}}) & =\left( {\frac{{y_{i}^{k}y_{j}^{k}}}{\left\langle {%
y^{k},y^{k}}\right\rangle }}\right) _{(i,j)\in\mathbb{N}\times\mathbb{N}%
}=\left( \left\langle P_{k}(e_{i}),e_{j}\right\rangle \right) _{(i,j)\in 
\mathbb{N}\times\mathbb{N}}, \\
A(\eta_{\mathbb{N}^{\ast}\setminus\left\{ k_{1},k_{2},\ldots,k_{n}\right\}
}) & =\left( \left\langle \left( Id-\sum_{s=1}^{n}P_{k_{s}}\right)
(e_{i}),e_{j}\right\rangle \right) _{(i,j)\in\mathbb{N}\times\mathbb{N}}.
\end{align*}

Using the above remark, we state the following theorem:

\begin{proposition}
There exists an isometric isomorphism between ${\mathcal{M}}$ and $\mathcal{S%
}_{Y}.$
\end{proposition}

\begin{theorem}
Each operator in $\mathcal{S}_{Y}$ is represented as a matrix whose entries
are integrals defined by scalar measures.
\end{theorem}

\section{The subalgebra $\mathcal{L}_{T}$}

In this section we study the smallest closed subalgebra with unity of $%
\mathcal{S}_{Y}\left( c_{0}\right) $ generated by a fixed element $%
T_{\lambda}\in$ $\mathcal{S}_{Y}\left( c_{0}\right) ,$ where $\lambda
=\left( \lambda_{n}\right) \in c_{0}$ and 
\begin{equation*}
T_{\lambda}=\sum_{n=1}^{\infty}\lambda_{n}P_{n}, 
\end{equation*}
we also show that, as $\mathcal{S}_{Y}\left( c_{0}\right) ,$ this algebra is
generated by a family of normal projections and, under certain conditions,
both algebras are isometrically isomorphic.

Let us denote by $\mathcal{L}_{T_{\lambda}}$ the closure of $alg_{\mathcal{S}%
_{Y}\left( c_{0}\right) }\left\{ Id,T_{\lambda}\right\} $ with respect to
the operator norm, that is, the closure of the space of polynomials in $%
T_{\lambda}$. Clearly, $\mathcal{L}_{T_{\lambda}}$ is a $C$-algebra since it
is closed Banach subalgebra of $\mathcal{S}_{Y}\left( c_{0}\right) $.

By Proposition \ref{P2}, $Sp\left( \mathcal{L}_{T_{\lambda}}\right) $ is compact
and, for each $H\in\mathcal{L}_{T_{\lambda}},$ where 
\begin{equation*}
\left\Vert H\right\Vert =\sup_{i\in\mathbb{N}}\left\Vert H\left(
e_{i}\right) \right\Vert =\left\Vert H\right\Vert _{sp}. 
\end{equation*}
On the other hand, since $\mathcal{S}_{Y}\left( c_{0}\right) $ has the power
multiplicative norm property, $\mathcal{L}_{T_{\lambda}}$ inherits such
property.

Under the conditions that $\mathcal{L}_{T_{\lambda}}$ is a $C$-algebra and $%
Sp\left( \mathcal{L}_{T_{\lambda}}\right) $ is compact, we conclude that $%
\mathcal{L}_{T_{\lambda}}$ is isometrically isomorphic to the space of all
continuous functions $C\left( Sp\left( \mathcal{L}_{T_{\lambda}}\right)
\right) $ provided by the supremum norm, that is, there exists an
isomorphism of algebras 
\begin{equation*}
\Psi:\mathcal{L}_{T_{\lambda}}\rightarrow C\left( Sp\left( \mathcal{L}%
_{T_{\lambda}}\right) \right) 
\end{equation*}
such that, for all $H\in\mathcal{L}_{T_{\lambda}},$ $\left\Vert H\right\Vert
=\left\Vert \Psi\left( H\right) \right\Vert _{\infty}.$

Suppose, for instant, that the range of the sequence $\lambda$ is infinite
and define the equivalence relation $n\sim
m\Leftrightarrow\lambda_{n}=\lambda _{m}.$ Observe that each equivalence
class of a non-null entry is at most a finite set. Let us denote by $\left\{
\lambda_{n_{1}},\lambda_{n_{2}},\ldots\right\} $ the collection of all
non-null representative of such classes. Of course, if $n_{1}<n_{2}<\ldots,$
then $\lim_{n_{i}\rightarrow \infty}\lambda_{n_{i}}=0=\lambda_{n_{0}}.$
Therefore, $\left\{ \lambda_{n_{1}},\lambda_{n_{2}},\ldots\right\} \cup\left\{
\lambda_{n_{0}}\right\} =\sigma\left( T_{\lambda}\right)$, the sequence of
all eigenvalues of $T_{\lambda}.$

Let us consider the unique homomorphism of algebra 
\begin{equation*}
\phi_{i}:alg_{\mathcal{S}_{Y}\left( c_{0}\right) }\left\{ Id,T_{\lambda
}\right\} \rightarrow\mathbb{K},\ \ i\in\mathbb{N}\cup\left\{ 0\right\} 
\end{equation*}
such that%
\begin{equation*}
\phi_{i}\left( T_{\lambda}\right) =\lambda_{n_{i}}. 
\end{equation*}
Thus, for any $H=\sum_{m=0}^{k}\alpha_{m}T_{\lambda}^{m}\in alg_{\mathcal{S}%
_{Y}\left( c_{0}\right) }\left\{ Id,T_{\lambda}\right\} ,$ we have%
\begin{align*}
\left\vert \phi_{i}\left( H\right) \right\vert & =\left\vert \alpha
_{0}+\sum_{m=1}^{k}\alpha_{m}\lambda_{n_{i}}^{m}\right\vert \leq\max\left\{
\left\vert \alpha_{0}\right\vert ,\left\vert
\sum_{m=1}^{k}\alpha_{m}\lambda_{n_{i}}^{m}\right\vert \right\} \\
& \leq\max\left\{ \left\vert \alpha_{0}\right\vert ,\left\Vert \sum
_{m=1}^{k}\alpha_{m}\lambda^{m}\right\Vert _{c_{0}}\right\} =\left\Vert
\left( \alpha_{0},\sum_{m=1}^{k}\alpha_{m}\lambda^{m}\right) \right\Vert
_{c_{0}^{+}} \\
& =\left\Vert
\alpha_{0}Id+T_{\sum_{m=1}^{k}\alpha_{m}\lambda^{m}}\right\Vert =\left\Vert
\sum_{m=0}^{k}\alpha_{m}T_{\lambda}^{m}\right\Vert =\left\Vert H\right\Vert ,
\end{align*}
that is, $\phi_{i}$ is continuous in $alg_{\mathcal{S}_{Y}\left(
c_{0}\right) }\left\{ Id,T_{\lambda}\right\} $. From this, $\phi_{i}$ can be
uniquely extended to a continuous homomorphism of algebra from $\mathcal{L}%
_{T_{\lambda}}$ into $\mathbb{K}.$

We claim that $\sigma\left( T_{\lambda}\right) ,$ equipped with a certain
topology, is homeomorphic to $Sp\left( \mathcal{L}_{T_{\lambda}}\right) .$

Note that the function%
\begin{equation*}
\Gamma:\sigma\left( T_{\lambda}\right) \rightarrow Sp\left( \mathcal{L}%
_{T_{\lambda}}\right) ;\ \lambda_{n_{i}}\longmapsto\Gamma\left(
\lambda_{n_{i}}\right) =\phi_{i}
\end{equation*}
is well-defined and is injective.

The next proposition tell us when an element of $\mathcal{S}_{Y}\left(
c_{0}\right) $ admits an inverse:

\begin{proposition}
Let $T\in\mathfrak{T}_{Y}(c_{0}).$ If $z\notin\sigma\left( T\right) ,$ then $%
zId-T$ is invertible in $\mathcal{S}_{Y}\left( c_{0}\right) .$

\begin{proof}
For $y\in R\left( zId-T\right) ,$ there exists $x\in c_{0}$ such that%
\begin{equation*}
\left( zId-T\right) \left( x\right) =y 
\end{equation*}
Since $z\notin\sigma\left( T\right) ,$ we can solve the above equation for $x
$ and get%
\begin{equation}\label{3.1}
x=\frac{1}{z}y+\frac{1}{z}Tx=\frac{1}{z}y+\frac{1}{z}\sum_{i=1}^{\infty
}\lambda_{i}\frac{\left\langle x,y^{\left( i\right) }\right\rangle }{%
\left\langle y^{\left( i\right) },y^{\left( i\right) }\right\rangle }%
y^{\left( i\right) }
\end{equation}
Applying the continuous functional $\left\langle \cdot,y^{\left( k\right)
}\right\rangle $ to $x$, we have%
\begin{align*}
\left\langle x,y^{\left( k\right) }\right\rangle & =\left\langle \frac {1}{z}%
y+\frac{1}{z}\sum_{i=1}^{\infty}\lambda_{i}\frac{\left\langle x,y^{\left(
i\right) }\right\rangle }{\left\langle y^{\left( i\right) },y^{\left(
i\right) }\right\rangle }y^{\left( i\right) },y^{\left( k\right)
}\right\rangle \\
& =\frac{1}{z}\left\langle y,y^{\left( k\right) }\right\rangle +\frac{1}{z}%
\lambda_{k}\left\langle x,y^{\left( k\right) }\right\rangle
\end{align*}
Now, solving the last equation for $\left\langle x,y^{\left( k\right)
}\right\rangle ,$ we obtain%
\begin{align*}
\left( 1-\frac{\lambda_{k}}{z}\right) \left\langle x,y^{\left( k\right)
}\right\rangle & =\frac{1}{z}\left\langle y,y^{\left( k\right) }\right\rangle
\\
\left\langle x,y^{\left( k\right) }\right\rangle & =\frac{1}{z-\lambda _{k}}%
\left\langle y,y^{\left( k\right) }\right\rangle
\end{align*}
Note that the sequence%
\begin{equation*}
\left( \frac{\lambda_{k}}{z-\lambda_{k}}\right) _{k\in\mathbb{N}}
\end{equation*}
is an element of $c_{0}.$ In fact, for a given $0<\epsilon<1$, there exists $%
i_{0}\in\mathbb{N}\,$\ such that%
\begin{equation*}
i\geq i_{0}\Rightarrow\left\vert \lambda_{i}\right\vert <\epsilon\left\vert
z\right\vert . 
\end{equation*}
Thus,%
\begin{equation*}
i\geq i_{0}\Rightarrow\left\vert \frac{\lambda_{i}}{z-\lambda_{i}}%
\right\vert =\frac{\left\vert \lambda_{i}\right\vert }{\left\vert
z\right\vert }<\epsilon. 
\end{equation*}
Now, replacing in \eqref{3.1}, we get%
\begin{align*}
x & =\frac{1}{z}y+\frac{1}{z}\sum_{i=1}^{\infty}\frac{\lambda_{i}}{%
z-\lambda_{i}}\frac{\left\langle y,y^{\left( i\right) }\right\rangle }{%
\left\langle y^{\left( i\right) },y^{\left( i\right) }\right\rangle }%
y^{\left( i\right) } \\
& =\frac{1}{z}y+\frac{1}{z}\sum_{i=1}^{\infty}\frac{\lambda_{i}}{%
z-\lambda_{i}}P_{i}\left( y\right) .
\end{align*}
Although $y$ belongs to $R\left( zId-T\right) ,$ the last expression holds
for any $y\in c_{0}.$ Thus, if we denote by 
\begin{equation*}
R_{z}\left( T\right) \left( y\right) =\frac{1}{z}y+\frac{1}{z}\sum
_{i=1}^{\infty}\frac{\lambda_{i}}{z-\lambda_{i}}P_{i}\left( y\right) , 
\end{equation*}
then $R_{z}\left( T\right) \left( \cdot\right) \in\mathcal{S}_{Y}\left(
c_{0}\right) ,$ since $\sum_{i=1}^{\infty}\frac{\lambda_{i}}{z-\lambda_{i}}%
P_{i}\left( \cdot\right) $ is compact and self-adjoint operator.\newline
Let us show that, effectively, $R_{z}\left( T\right) \left( \cdot\right) $
is the inverse operator of $zId-T:$%
\begin{align*}
& \left. \left[ \left( zId-T\right) \circ R_{z}\left( T\right) \right]
\left( y\right) =\right. \left( zId-T\right) \left( \frac{1}{z}y+\frac{1}{z}%
\sum_{i=1}^{\infty}\frac{\lambda_{i}}{z-\lambda_{i}}P_{i}\left( y\right)
\right) \\
& =y+\sum_{i=1}^{\infty}\frac{\lambda_{i}}{z-\lambda_{i}}P_{i}\left(
y\right) -\frac{1}{z}\sum_{i=1}^{\infty}\lambda_{i}P_{i}\left( y\right) -%
\frac{1}{z}\sum_{i=1}^{\infty}\frac{\lambda_{i}^{2}}{z-\lambda_{i}}%
P_{i}\left( y\right) ;\ \ T\left( P_{i}\left( y\right) \right)
=\lambda_{i}P_{i}\left( y\right) \\
& =y+\sum_{i=1}^{\infty}\left[ \underset{=0}{\underbrace{\frac{\lambda_{i}}{%
z-\lambda_{i}}-\frac{\lambda_{i}}{z}-\frac{\lambda_{i}^{2}}{z\left(
z-\lambda_{i}\right) }}}\right] P_{i}\left( y\right) =y=Id\left( y\right)
\end{align*}
In the other direction, since 
\begin{equation*}
P_{j}\circ P_{i}\left( x\right) =\left\{ 
\begin{array}{ccc}
P_{i}\left( x\right) & if & j=i \\ 
0 & if & j\neq i%
\end{array}
\right. , 
\end{equation*}
we have 
\begin{align*}
& \left. \left[ R_{z}\left( T\right) \circ\left( zId-T\right) \right] \left(
x\right) =\right. zR_{z}\left( T\right) \left( x\right) -R_{z}\left(
T\right) \left( Tx\right) \\
& =x+\sum_{i=1}^{\infty}\frac{\lambda_{i}}{z-\lambda_{i}}P_{i}\left(
x\right) -\sum_{i=1}^{\infty}\lambda_{i}R_{z}\left( T\right) \left(
P_{i}\left( x\right) \right) \\
& =x+\sum_{i=1}^{\infty}\frac{\lambda_{i}}{z-\lambda_{i}}P_{i}\left(
x\right) -\sum_{i=1}^{\infty}\lambda_{i}\left[ \frac{1}{z}P_{i}\left(
x\right) +\frac{1}{z}\sum_{j=1}^{\infty}\frac{\lambda_{j}}{z-\lambda_{j}}%
P_{j}\left( P_{i}\left( x\right) \right) \right] \\
& =x+\sum_{i=1}^{\infty}\left[ \underset{=0}{\underbrace{\frac{\lambda_{i}}{%
z-\lambda_{i}}-\frac{\lambda_{i}}{z}-\frac{\lambda_{i}^{2}}{z\left(
z-\lambda_{i}\right) }}}\right] P_{i}\left( x\right) =x=Id\left( x\right)
\end{align*}
Therefore, $R_{z}\left( T\right) =\left( zId-T\right) ^{-1}\in \mathcal{S}%
_{Y}\left( c_{0}\right) .$
\end{proof}
\end{proposition}

\begin{corollary}
If $z\notin\sigma\left( T_{\lambda}\right) ,$ then $R_{z}\left( T_{\lambda
}\right) =\left( zId-T_{\lambda}\right) ^{-1}\in\mathcal{L}_{T_{\lambda}}$.

\begin{proof}
We already know that $\mathcal{S}_{Y}\left( c_{0}\right) $ is a $C$-algebra
with unity and $zId-T_{\lambda }$ is invertible in $\mathcal{S}_{Y}\left(
c_{0}\right) .$ By Th. 6.10 in \cite{RO}, we have that 
\begin{equation*}
R_{z}\left( T_{\lambda }\right) \in \overline{alg_{\mathcal{S}_{Y}\left(
c_{0}\right) }}\left\{ Id,zId-T_{\lambda }\right\} 
\end{equation*}%
Now, since $\overline{alg_{\mathcal{S}_{Y}\left( c_{0}\right) }}\left\{
Id,zId-T_{\lambda }\right\} $ is the smallest closed subalgebra that
contains to $zId-T_{\lambda }$, we conclude that $R_{z}\left( T_{\lambda
}\right) \in \mathcal{L}_{T_{\lambda }}.$
\end{proof}
\end{corollary}

\begin{proposition}
The function $\Gamma$ is bijective.

\begin{proof}
By above, $\Gamma$ is injective. If $\phi\in Sp\left( \mathcal{L}%
_{T_{\lambda}}\right) ,$ then $\phi\left( T_{\lambda}\right) =z,$ for some $%
z\in\mathbb{K}.$ Suppose that $z\notin\sigma\left( T_{\lambda}\right) ,$
hence $zId-T_{\lambda}$ has an inverse and, by the previous corollary, $%
R_{z}\left( T_{\lambda}\right) \in\mathcal{L}_{T_{\lambda}}.$ Since the
function $\phi$ is a homomorphism between algebras with unities, we have%
\begin{equation*}
1=\phi\left( Id\right) =\phi\left( \left( zId-T_{\lambda}\right)
^{-1}\circ\left( zId-T_{\lambda}\right) \right) =\phi\left( \left(
zId-T_{\lambda}\right) ^{-1}\right) \phi\left( zId-T_{\lambda}\right) , 
\end{equation*}
but, by the linearity of $\phi,$ the factor $\phi\left( zId-T_{\lambda
}\right) $ is null$,$ which is a contradiction. Thus, if $\phi\in Sp\left( 
\mathcal{L}_{T_{\lambda}}\right) ,$ then there exists $\mu\in\sigma\left(
T_{\lambda}\right) $ such that $\phi=\phi_{\mu}$ and therefore $\Gamma$ is
bijective.
\end{proof}
\end{proposition}

\begin{remark}
We have identified $Sp\left( \mathcal{L}_{T_{\lambda}}\right) $ with $%
\sigma\left( T_{\lambda}\right) $ through the bijective function $\Gamma$.
Let us consider the induced topology by $\mathbb{K}$ on $\sigma\left(
T_{\lambda}\right) .$ Note that $\sigma\left( T_{\lambda}\right) $ is
compact.
\end{remark}

\begin{proposition}
$\sigma\left( T_{\lambda}\right) $ is homeomorphic to $Sp\left( \mathcal{L}%
_{T_{\lambda}}\right) $

\begin{proof}
We claim first that $\Upsilon=\Gamma^{-1}$ is continuous. In fact, if $%
\phi_{\alpha}\rightarrow\phi$ in the induced topology on $Sp\left( \mathcal{L%
}_{T_{\lambda}}\right) $ by the product topology in $\mathbb{K}^{\mathcal{L}%
_{T_{\lambda}}},$ then%
\begin{equation*}
\phi_{\alpha}\left( H\right) \rightarrow\phi\left( H\right) 
\end{equation*}
for each $H\in\mathcal{L}_{T_{\lambda}}.$ In particular, 
\begin{equation*}
\phi_{\alpha}\left( T_{\lambda}\right) \rightarrow\phi\left( T_{\lambda
}\right) 
\end{equation*}
or, equivalently, 
\begin{equation*}
\Upsilon\left( \phi_{\alpha}\right) \rightarrow\Upsilon\left( \phi\right) . 
\end{equation*}
Now, since $\Upsilon$ is bijective and continuous, $Sp\left( \mathcal{L}%
_{T_{\lambda}}\right) $ is compact and $\sigma\left( T_{\lambda}\right) $ is
a Hausdorff space, we conclude that $\Upsilon$ is a homeomorphism.\newline
\end{proof}
\end{proposition}

By this proposition and by the uniqueness of $X$ (up to homeomorphism) for
which $\mathcal{L}_{T_{\lambda}}\cong C\left( X\right) ,$ we have 
\begin{equation*}
\mathcal{L}_{T_{\lambda}}\cong C\left( Sp\left( \mathcal{L}%
_{T_{\lambda}}\right) \right) \cong C\left( \sigma\left( T_{\lambda}\right)
\right) . 
\end{equation*}

Let us identify the isometric isomorphism $\Psi.$ Replacing $Sp\left( 
\mathcal{L}_{T_{\lambda}}\right) $ by $\sigma\left( T_{\lambda}\right) $ and
considering the map 
\begin{equation*}
G_{T_{\lambda}}\overset{notation}{=}f_{T_{\lambda}}:\sigma\left( T_{\lambda
}\right) \rightarrow\mathbb{K};\ \lambda_{n_{i}}\longmapsto
f_{T_{\lambda}}\left( \lambda_{n_{i}}\right) =\lambda_{n_{i}},\ 
\end{equation*}
we can get, for a fixed $H=\alpha_{0}Id+\sum_{m=1}^{k}\alpha_{m}T_{\lambda
}^{m}\in alg_{\mathcal{S}_{Y}\left( c_{0}\right) }\left\{ Id,T_{\lambda
}\right\} ,$ the map $G_{H}:\sigma\left( T_{\lambda}\right) \rightarrow\mathbb{K}$ defined by
\begin{equation*}
G_{H}\left( \lambda_{n_{i}}\right) =\left\{ 
\begin{array}{ccc}
\alpha_{0}+\sum_{m=1}^{k}\alpha_{m}\lambda_{n_{i}}^{m}=\alpha_{0}+\sum
_{m=1}^{k}\alpha_{m}\left[ f_{T_{\lambda}}\left( \lambda_{n_{i}}\right) %
\right] ^{m} & if & i\in\mathbb{N} \\ 
\alpha_{0} & if & i=0%
\end{array}
\right. . 
\end{equation*}
From this, we can define 
\begin{equation*}
G:alg_{\mathcal{S}_{Y}\left( c_{0}\right) }\left\{ Id,T_{\lambda}\right\}
\rightarrow C\left( \sigma\left( T_{\lambda}\right) \right) ;\ H\longmapsto
G_{H}, 
\end{equation*}
the well-known Gelfand transformation. Clearly, $G$ is a homomorphism of
algebras, $\left\Vert H\right\Vert _{sp}=\left\Vert G_{H}\right\Vert $ for
any $H\in alg_{\mathcal{S}_{Y}\left( c_{0}\right) }\left\{ Id,T_{\lambda
}\right\} ,$ and since $\mathcal{L}_{T_{\lambda}}$ satisfies the condition
given by Proposition \ref{P2}, we have%
\begin{equation*}
\left\Vert H\right\Vert _{sp}=\left\Vert H\right\Vert . 
\end{equation*}
Thus, $G$ is an isometry and then it can be extended to the whole $\mathcal{L%
}_{T_{\lambda}}.$ Let us denote by the same capital letter $G$ such
extension.

\begin{proposition}
$G$ is an isometric isomorphism of algebras.

\begin{proof}
It is enough to prove that $G$ is surjective. By the fact that $G$ is a
homomorphism of algebras and the image of $T_{\lambda }$ by $G$ is the
identity map $G_{T_{\lambda }}=f_{T_{\lambda }}$, the collection $\left\{
1,f_{T_{\lambda }},f_{T_{\lambda }}^{2},\ldots \right\} $ is the image of $%
\left\{ Id,T_{\lambda },T_{\lambda }^{2},\ldots \right\} $ by $G$. Now, by the compactness of $\sigma \left( T_{\lambda }\right) $, Theorem
5.28 in \cite{RO} guarantees that $alg_{C\left( \sigma \left(
T_{\lambda }\right) \right) }\left\{ 1,f_{T_{\lambda }}\right\} $ is dense
in $C\left( \sigma \left( T_{\lambda }\right) \right) .$ Thus, if $f\in
C\left( \sigma \left( T_{\lambda }\right) \right) ,$ then there exists a
sequence $\left\{ g_{n}\right\} _{n\in \mathbb{N}}$ in $alg_{C\left( \sigma
\left( T_{\lambda }\right) \right) }\left\{ 1,f_{T_{\lambda }}\right\} $
such that $f=\lim_{n\rightarrow \infty }g_{n}.$ Now, for each $n\in \mathbb{N%
},$ there exists $H_{n}\in \mathcal{L}_{T_{\lambda }}$ such that%
\begin{equation*}
G_{H_{n}}=g_{n}.
\end{equation*}%
By the fact that $G$ is an isometry$,$ the sequence $\left( H_{n}\right) $
is a Cauchy sequence in $\mathcal{L}_{T_{\lambda }}$ and then convergence to
an $H\in \mathcal{L}_{T_{\lambda }}$. Since $G$ is continuous, we have that 
\begin{equation*}
G_{H}=\lim_{n\rightarrow \infty }G_{H_{n}}=\lim_{n\rightarrow \infty
}g_{n}=f.
\end{equation*}
\end{proof}
\end{proposition}

\begin{remark}
By above proposition, the Gelfand transformation $G$ is an isometric
isomorphism between $\mathcal{L}_{T_{\lambda}}$ and $C\left( Sp\left( 
\mathcal{L}_{T_{\lambda}}\right) \right) $ or $C\left( \sigma\left(
T_{\lambda}\right) \right) .$
\end{remark}

\begin{remark}
If we suppose that, for $\lambda,\mu\in c_{0},$ the corresponding sets $%
\left\{ \lambda_{n}:n\in\mathbb{N}\right\} $ and $\left\{ \mu_{n}:n\in%
\mathbb{N}\right\} $ are infinite, then both are connected by a bijective
correspondence which is a homeomorphism if we provide them with the discrete
topology. Therefore, the spaces $\sigma\left( T_{\lambda}\right) ,$ $%
\sigma\left( T_{\mu}\right) $ and $\mathbb{N}^{\ast}$ are homeomorphic each
other. Now, since any C-algebra with unity is isometrically isomorphic to $%
C\left( X\right) ,$ where $X$ is compact space and it is unique up to
homeomorphism, we conclude that $\mathcal{L}_{T_{\lambda}}\cong\mathcal{L}%
_{T_{\mu}}\cong\mathcal{S}_{Y}.$
\end{remark}

\subsection{Spectral measure}

By the previous section, there exists an isometric isomorphism of algebras $%
\Phi=G^{-1}:C\left( \sigma\left( T_{\lambda}\right) \right) \rightarrow%
\mathcal{L}_{T_{\lambda}}.$ Let us denote by $\Omega\left( \sigma\left(
T_{\lambda}\right) \right) $ the Boolean ring of all clopen subsets of $%
\sigma\left( T_{\lambda}\right) .$ Of course, $\eta_{C}$ is continuous if
and only if $C\in\Omega\left( \sigma\left( T_{\lambda}\right) \right) .$

Note that the elements of $\Omega\left( \sigma\left( T_{\lambda}\right)
\right) $ can be classified as follows: the first type are those which are
finite subsets of $\sigma\left( T_{\lambda}\right) \setminus\left\{
\lambda_{n_{0}}\right\} $ and the second type are those which are complement
on $\sigma\left( T_{\lambda}\right) $ of the first type.

Now, since $\Phi$ is a homomorphism of algebras, we have 
\begin{equation*}
\Phi\left( \eta_{C}\right) =\Phi\left( \eta_{C}^{2}\right) =\Phi\left(
\eta_{C}\right) ^{2}. 
\end{equation*}
In other words, $\Phi\left( \eta_{C}\right) $ is a projection in $\mathcal{L}%
_{T_{\lambda}}$ and if $C\in\Omega\left( \sigma\left( T_{\lambda}\right)
\right) \setminus\left\{ \varnothing\right\} ,$ then $\Phi\left(
\eta_{C}\right) $ is a non-null$.$

On the other hand, by the fact that the linear hull of $\left\{
\eta_{C}:C\in\Omega\left( \sigma\left( T_{\lambda}\right) \right) \right\} $
is dense in $C\left( \sigma\left( T_{\lambda}\right) \right) $, for any
fixed $f\in C\left( \sigma\left( T_{\lambda}\right) \right) $ and $%
\epsilon>0,$ there exists a finite clopen partition $\left\{
C_{k}:k=1,\ldots,s\right\} $ of $\sigma\left( T_{\lambda}\right) $ and a
finite collection of scalars $\left\{ \alpha_{k}:k=1,\ldots,s\right\} $ such
that 
\begin{equation}\label{3.2}
\left\Vert f-\sum_{k=1}^{s}\alpha_{k}\eta_{C_{k}}\right\Vert
_{\infty}=\sup_{x\in\sigma\left( T_{\lambda}\right) }\left\vert f\left(
x\right) -\sum_{k=1}^{s}\alpha_{k}\eta_{C_{k}}\left( x\right) \right\vert
<\epsilon
\end{equation}
Since $\left\{ C_{k}:k=1,\ldots,s\right\} $ is a clopen partition of $%
\sigma\left( T_{\lambda}\right) ,$ only one of these sets is of second type,
say $C_{1}=\sigma\left( T_{\lambda}\right) \setminus\left\{
\lambda_{m_{1}},\lambda_{m_{2}},\ldots,\lambda_{m_{n}}\right\} .$ From this, 
$\cup_{k=2}^{s}C_{k}=\left\{ \lambda_{m_{1}},\lambda_{m_{2}},\ldots
,\lambda_{m_{n}}\right\} .$ Using the characteristic functions properties
and the fact that the single subsets $\left\{ \lambda_{k}\right\} $ belong
to $\Omega\left( \sigma\left( T_{\lambda}\right) \right) $, we can rewrite \eqref{3.2} as follows: 
\begin{align*}
& \left\Vert f-\left[ \alpha_{1}\eta_{\sigma\left( T_{\lambda}\right)
}+\sum_{l=1}^{n}\left( \alpha_{l}-\alpha_{1}\right) \eta_{\left\{
\lambda_{m_{l}}\right\} }\right] \right\Vert _{\infty} \\
& =\sup_{x\in\sigma\left( T_{\lambda}\right) }\left\vert f\left( x\right) -
\left[ \alpha_{1}\eta_{\sigma\left( T_{\lambda}\right) }\left( x\right)
+\sum_{l=1}^{n}\left( \alpha_{l}-\alpha_{1}\right) \eta_{\left\{
\lambda_{m_{l}}\right\} }\left( x\right) \right] \right\vert <\epsilon
\end{align*}
and without loss of generality, we can assume that 
\begin{align*}
& \left\Vert f-\left[ \alpha_{1}\eta_{\sigma\left( T_{\lambda}\right)
}+\sum_{l=1}^{n}\left( \alpha_{l}-\alpha_{1}\right) \eta_{\left\{
\lambda_{m_{l}}\right\} }\right] \right\Vert _{\infty} \\
& =\sup_{x\in\sigma\left( T_{\lambda}\right) }\left\vert f\left( x\right) -
\left[ f\left( \lambda_{n_{0}}\right) \eta_{\sigma\left( T_{\lambda }\right)
}\left( x\right) +\sum_{l=1}^{n}\left[ f\left( \lambda_{m_{l}}\right)
-f\left( \lambda_{n_{0}}\right) \right] \eta_{\left\{
\lambda_{m_{l}}\right\} }\left( x\right) \right] \right\vert <\epsilon
\end{align*}
Using the isometry $\Phi,$ we have
\begin{equation}\label{3.3}
\left\Vert \Phi\left( f\right) -\left[ f\left( \lambda_{n_{0}}\right)
Id+\sum_{l=1}^{n}\left[ f\left( \lambda_{m_{l}}\right) -f\left(
\lambda_{n_{0}}\right) \right] E_{l}\right] \right\Vert <\epsilon,
\end{equation}
where $E_{l}$ is the corresponding projection $\Phi\left( \eta_{\left\{
\lambda_{m_{l}}\right\} }\right) .$ At the same time, \eqref{3.3} shows that the space generated by $\left\{ E\in\mathcal{L}_{T_{\lambda}}:E^{2}=E\right\}$ is dense in $\mathcal{L}_{T_{\lambda}}.$

Let us consider the following set-function:%
\begin{equation*}
m_{T_{\lambda}}:\Omega\left( \sigma\left( T_{\lambda}\right) \right)
\rightarrow\mathcal{L}_{T_{\lambda}};\ \ C\longmapsto m_{T_{\lambda}}\left(
C\right) =\Phi\left( \eta_{C}\right) =E_{C}. 
\end{equation*}
In similar way as in subsection \ref{S2.2}, $m_{T_{\lambda}}$ is a finite additive
measure valued-projection which is known as spectral measure associated to $%
\mathcal{L}_{T_{\lambda}}$.

Since $E$ is, in particular, an element of $\mathcal{S}_{Y},$ there exists $%
\alpha\in\mathbb{K}$ and $\mu=\left( \mu_{i}\right) \in c_{0}$ such that $%
E=\alpha Id+T_{\mu}.$ By the fact that $E$ is a projection, $E^{2}=E,$ that
is, $\left( \alpha Id+T_{\mu}\right) \left( \alpha Id+T_{\mu}\right) =\alpha
Id+T_{\mu}$ or equivalent to $\left( \alpha^{2}-\alpha\right) Id+T_{\left(
2\alpha-1\right) \mu+\mu^{2}}=0.$\ Taking the norm of this operator, we get%
\begin{equation*}
0=\left\Vert \left( \alpha^{2}-\alpha\right) Id+T_{\left( 2\alpha-1\right)
\mu+\mu^{2}}\right\Vert =\max\left\{ \left\vert \alpha\left( \alpha
-1\right) \right\vert ,\left\Vert \left( 2\alpha-1\right) \mu+\mu
^{2}\right\Vert \right\} 
\end{equation*}
From this, if $\alpha=0,$ then $\mu_{i}=0$ or $\mu_{i}=1.$ Since $%
\lim_{n\rightarrow\infty}\mu_{n}=0,$ $\mu_{i}=0$ for all $i\in\mathbb{N},$
excepts for a finite collection $\left\{ i_{1},i_{2},\ldots,i_{n}\right\} $
for which $\mu_{i_{s}}=1,\ s=1,2,\ldots,n.$ Thus, $E=\sum_{s=1}^{n}P_{i_{s}}.
$\newline
On the other hand, if $\alpha=1,$ then we also get finite collection, say $%
\left\{ j_{1},j_{2},\ldots,j_{m}\right\} $ for which $\mu_{j_{k}}=-1,\
k=1,2,\ldots,m,$ and the rest of the elements of the sequence $\mu$ are $0.$
Thus, $E=Id-\sum_{k=1}^{m}P_{j_{k}}.$\newline
If $\sum_{s=1}^{n}P_{s}\in\mathcal{L}_{T_{\lambda}},\ n\neq1,$ then $%
P_{s}\notin\mathcal{L}_{T_{\lambda}},$ for any $s\in\left\{
1,2,\ldots,n\right\} .$

\begin{remark}
Note that the length of the sum in $E=\sum_{s=1}^{n}P_{s}$ depends on
exclusively for the repetition of some non-null entries of the sequence $%
\lambda.$
\end{remark}

Following the same arguments as in subsection \ref{S2.2}, for $f\in C\left(
\sigma\left( T_{\lambda}\right) \right) $ and $\alpha=\left\{
C_{1},C_{2},\ldots,C_{n};x_{1},x_{2},\ldots,x_{n}\right\} \in\mathcal{D},$
where $\sigma\left( T_{\lambda}\right) =\sqcup_{k=1}^{n}C_{k},$ we define%
\begin{equation*}
\omega_{\alpha}(f,m_{T_{\lambda}},\sigma\left( T_{\lambda}\right)
)=\sum_{k=1}^{n}f\left( x_{k}\right) m_{T_{\lambda}}\left( C_{k}\right)
=\sum_{k=1}^{n}f\left( x_{k}\right) E_{C_{k}}. 
\end{equation*}
and since the function $f$ can be reached by a net 
\begin{equation*}
\left\{ \sum_{k=1}^{n}f\left( x_{k}\right) \eta_{C_{k}}\right\} _{\alpha\in%
\mathcal{D}}
\end{equation*}
in $C\left( \sigma\left( T_{\lambda}\right) \right) ,$ the isometry of $\Phi$
allows us to get%
\begin{equation*}
\underset{\alpha\in\mathcal{D}}{\lim}\omega_{\alpha}(f,m_{T_{\lambda}},%
\sigma\left( T_{\lambda}\right) )=\Phi\left( f\right) 
\end{equation*}
Therefore, the operator $\Phi\left( f\right) $ is interpreted as an
integral, that is, 
\begin{equation*}
\Phi\left( f\right) =\int_{\sigma\left( T_{\lambda}\right)
}fdm_{T_{\lambda}}=\lim_{\alpha\in\mathcal{D}}\omega_{\alpha}(f,m_{T_{%
\lambda}},\sigma\left( T_{\lambda}\right) ) 
\end{equation*}
For example, for $f_{T_{\lambda}},$ $\eta_{\left\{ \lambda_{n}\right\} }$ or 
$f\equiv1,$ their respective integral are 
\begin{equation*}
T_{\lambda}=\Phi\left( f_{T_{\lambda}}\right) =\int_{\sigma\left(
T_{\lambda}\right) }f_{T_{\lambda}}dm_{T_{\lambda}};\ \
E_{n}=\int_{\sigma\left( T_{\lambda}\right) }\eta_{\left\{
\lambda_{n}\right\} }dm_{T_{\lambda}}=\left\{ 
\begin{array}{c}
\sum_{s=1}^{n}P_{s} \\ 
or \\ 
Id-\sum_{s=1}^{m}P_{s}%
\end{array}
\right.\
\end{equation*}
$$Id=\Phi\left( 1\right) =\int_{\sigma\left(
T_{\lambda}\right) }dm_{T_{\lambda}}.$$

\begin{remark}
Since $\left\{ P_{n}:n\in\mathbb{N}\right\} $ is a family of normal
projections, we conclude that $\left\{ E_{k}:k\in\mathbb{N}\right\} $ is
also a family of normal projections. Even more, using the inner product for
operators in $\mathcal{A}_{1},$ $\left\{ E_{k}:k\in\mathbb{N}\right\} $ is
an orthonormal family.\newline
Now, since $\left\{ \eta_{\left\{ \lambda _{k}\right\} }:k\in\mathbb{N}%
\right\} $ generates to $C\left( \sigma\left( T_{\lambda}\right) \right) ,$
the isometry isomorphism $\Phi$ tells us $\mathcal{L}_{T_{\lambda}}$ is
generated by \thinspace$\left\{ Id, E_1, E_2, \ldots\right\}$. In other
words, $\mathcal{L}_{T_{\lambda}}=\overline{Span\left\{
Id,E_{1},E_{2},\ldots\right\} }.$\newline
According to this, $\mathcal{L}_{T_{\lambda}}$ has the same structure than $%
\mathcal{S}_{Y}$ and, therefore, with the same arguments developed in
subsection \ref{S2.3}, each operator in $\mathcal{L}_{T_{\lambda}}$ admits a matrix
representation whose entries are integrals defined by scalar measures.
\end{remark}

\textbf{Acknowledgement.} We would like to thank Professor Bertin Diarra for
his useful suggestions and remarks.

\bibliographystyle{amsplain}

\bigskip\bigskip

\end{document}